# On the period of Pell-Narayana sequence in some groups

Bahar Kuloğlu [1], Engin Özkan[2], Marin Marin[3,*]

[1] *Department of Mathematics, Graduate School of Natural and Applied Sciences, Erzincan Binali Yıldırım University, Yalnızbağ Campus, 24100, Erzincan, Turkey.*
[2] *Department of Mathematics, Faculty of Arts and Sciences, Erzincan Binali Yıldırım University, Yalnızbağ Campus, 24100, Erzincan, Turkey*
[3] *Department of Mathematics and Computer Science Transilvania University of Brasov, 500036 Brasov Romania*

[*]*Corresponding author:* m.marin@unitbv.ro (M. Marin)

**ABSTRACT.** In this paper, the Pell-Narayana sequence modulo $m$ is studied. The paper outlines the definition of Pell-Narayana numbers and some of their combinatorial links with Eulerian, Catalan and Delannoy numbers and other special functions. From the definition, the Pell-Narayana orbit of a 2-generator group for a generating pair $(x, y) \in G$ is defined, so that the lengths of the period of the Pell-Narayana orbit can be examined. These yield in turn the Pell-Narayana lengths of the polyhedral group and the binary polyhedral group for the generating pair $(x, y)$ and associated properties. Also, the period of the Pell-Narayana orbit of the groups $Q_8, Q_8 \times \mathbb{Z}_{2m}$ and $Q_8 \times_\varphi \mathbb{Z}_{2m}, for\ m \geq 3$ were obtained.

**Keywords:** Pell-Narayana sequence, Pell-Narayana length, Pell-Narayana orbits, polyhedral group, binary polyhedral group.

**MSC2020:** 11B39, 11B50, 11B75.

## 1. INTRODUCTION

Since Fibonacci numbers exist in nature and find application in many branches of science, they have managed to attract the attention of many scientists for many years. Inspired by the Fibonacci numbers, many number sequences were defined and their properties were examined. One of the most important of these number sequences is the Narayana numbers. Narayana numbers are also associated with many other number sequences and their properties are studied [4,16-18].

The Pell-Narayana sequence is given by [1]:

$$PN_{n+3} = 2PN_{n+2} + PN_n, \qquad (1.1)$$

with the initial values $PN_0 = 0, PN_1 = 1$ and $PN_2 = 1$. That is, the Pell-Narayana sequence is $\{0,1,1,2,5,11,24, \dots\}$ (A078012) [22].

Cyclic group applications on Fibonacci number sequences, which started with Wall [25], continued to be studied with different recurrences on different numbers [6,8,12-14]. Later in [7,11], these applications were applied to 3-step Fibonacci numbers (Tribonacci sequence) and in [10] these applications were moved to $k$-step Fibonacci numbers.

Simple groups of order less than a million were considered by Campbell, Doostie and Robertson [2] and the binary polyhedral and binary polyhedral groups were studied in [24] and in [3]. In addition to these studies, Knox, Kumari, Aküzüm [9,20,21] expanded the study to the $k$-nacci sequences in finite groups. Similarly, it is possible to encounter the studies of these studies on 3-recurrences sequences such as Pell, Pell-Padovan, Jacobsthal-Padovan with their extended concepts on finite groups and fields [5,27,28,15].

It is aimed to examine the relationship between number sequences and algebra more closely by expanding the boundaries of working on number sequences in general, identities showing the relations related to sequences, Binet formula, and generating function.

Some similar results can be found in [29-30]. In this study, the connection of the Pell-Narayana number sequence with algebra is tried to be discussed by starting from the defined relation and groups.

We next extend the concept to the Pell-Narayana sequences.

In section 2, we examine Pell-Narayana sequence in finite groups and researched the period of the Pell-Narayana sequence according to modulo $m$. In section 3, we study in Pell-Narayana sequence orbits of finite groups. Finally, we examine Pell-Narayana sequence lengths of generating pairs in groups and we made applications to better understand these studies.

## 2. The Pell-Narayana sequences

A Pell-Narayana sequence $\{PN_r\}$ is defined [1] recursively by the equation

$$PN_r = 2PN_{r-1} + PN_{r-3}, \tag{1}$$

with the initial condition $PN_0 = 0, PN_1 = 1, PN_2 = 1$.

Kalman [26] pointed out that all sequences with 2, 3, ... recurrences are essentially a linear combination of the preceding $k$ terms of the equation defined below:

$$a_{n+k} = c_0 a_n + c_1 a_{n+1} + \cdots + c_{k-1} a_{n+k-1},$$

$\{c_0, c_1, \ldots, c_{k-1}\} \in \mathbb{R}$. In [26], these coefficients with companion matrix can be represented in closed form as follows:

$$A_k = [a_{ij}]_{k \times k} = \begin{bmatrix} 0 & 1 & 0 & \cdots & 0 & 0 \\ 0 & 0 & 1 & \cdots & 0 & 0 \\ 0 & 0 & 0 & \cdots & 0 & 0 \\ \vdots & \vdots & \vdots & \cdots & \vdots & \vdots \\ 0 & 0 & 0 & \cdots & 0 & 1 \\ c_0 & c_1 & c_2 & \cdots & c_{k-2} & c_{k-1} \end{bmatrix}_{k \times k}. \quad (2)$$

Then by an inductively,

$$A_k^n = \begin{bmatrix} a_0 \\ a_1 \\ \vdots \\ a_{k-1} \end{bmatrix} = \begin{bmatrix} a_n \\ a_{n+1} \\ \vdots \\ a_{n+k-1} \end{bmatrix}. \quad (3)$$

If the first $k$ terms of a sequence are repeated in the sequence from the beginning, the sequence is called simply the periodic sequence, and the period of this sequence is also $k$. For example, the sequence, $x_1, x_2, x_3, x_4, x_5, x_1, x_2, x_3, x_4, x_5 \ldots$ , is simply periodic with period 5. On the other hand, if consecutive $k$ terms repeat in the sequence between terms of a sequence, this sequence is called a periodic sequence, and the period of this sequence is $k$.

**Definition 1.** Consider a finite generated group $G = \langle A \rangle$ for the finite set $A = \{a_1, a_2, \ldots, a_n\}$. For $0 \leq i \leq n-1$, if $x_i = a_{i+1}$ and for $i \geq 0$ if $x_{i+n} = \prod_{j=1}^{n} x_{i+j-1}$ then the sequence $\{x_0, x_1, \ldots\}$ is called the Fibonacci orbit of $G$ and denoted by $F_A(G)$. If $F_A(G)$ is a simply periodic sequence, then the period of the sequence is called the Fibonacci length of $G$, shown $LEN_A(G)$ [3].

## 2.1 The modulo $m$

From statements (1) and (3), we can write

$$\begin{bmatrix} PN_{n-2} \\ PN_{n-1} \\ PN_n \end{bmatrix} = \begin{bmatrix} 0 & 1 & 0 \\ 0 & 0 & 1 \\ 1 & 0 & 2 \end{bmatrix} \begin{bmatrix} PN_{n-3} \\ PN_{n-2} \\ PN_{n-1} \end{bmatrix}, \quad (4)$$

for the Pell-Narayana sequence. Let us take

$$K = [k_{ij}]_{3 \times 3} = \begin{bmatrix} 0 & 1 & 0 \\ 0 & 0 & 1 \\ 1 & 0 & 2 \end{bmatrix}.$$

which is said to be Pell-Narayana matrix. Then by an inductively, for $n \geq 0$,

$$K^n \begin{bmatrix} 0 \\ 1 \\ 1 \end{bmatrix} = \begin{bmatrix} PN_n \\ PN_{n+1} \\ PN_{n+2} \end{bmatrix}, \quad (5)$$

is obtained.

Now, let's reduce each term of the Pell-Narayana sequence to mod $m$ and denote the resulting new sequence with $\{PN_n^{(m)}\}$. That is,

$$\{PN_n^{(m)}\} = \{PN_0^{(m)}, PN_1^{(m)}, PN_2^{(m)}, \ldots\}.$$

which has the same recurrence relation as in (1). These recurrences belong to a family of third order lacunary-type Padovan, Perrin, Pell-Narayana, Narayana and Plastic sequences which have been extensively explored by Anderson, Horadam, Petroudi, Soykan and Shannon [1,19,23].

**Theorem 2.** $\{PN_n^{(m)}\}$ forms a simply periodic sequence.

**Proof.** Since there are only a finite number $m^3$ of possible term triplets, the sequence repeats, and repeating the triple results in iteration of all subsequent terms.

From (1), we have

$$PN_{n+3} = 2PN_{n+2} + PN_n,$$

so, if

$$PN_{i+2}^{(m)} = PN_{j+2}^{(m)},$$
$$PN_{i+1}^{(m)} = PN_{j+1}^{(m)},$$
$$PN_i^{(m)} = PN_j^{(m)},$$

then

$$PN_{i-j+2}^{(m)} = PN_2^{(m)}, PN_{i-j+1}^{(m)} = PN_1^{(m)} \text{ and } PN_{i-j}^{(m)} = PN_0^{(m)},$$

which implies that the sequence $\{PN_n^{(m)}\}$ is a simply periodic, as required.

Now, let $\text{Per}(PN^{(m)})$ show the smallest period of the sequence $\{PN_n^{(m)}\}$.

Let $p_i$ be distinct primes. If $m = \prod_{i=1}^{t} p_i^{e_i} (t \geq 1)$ then we get

$\text{Per}(PN^{(m)}) = \text{lcm}[\text{Per}(PN^{(p_i^{e_i})})]$, the least common multiple of the $\text{Per}(PN^{(p_i^{e_i})})$.

**Theorem 3.** If $m = \prod_{i=1}^{t} p_i^{e_i} (t \geq 1)$ where the $p_i$ are distinct primes, then $\text{Per}(PN^{(m)}) = \text{lcm}[\text{Per}(PN^{(p_i^{e_i})})]$, the least common multiple of the $\text{Per}(PN^{(p_i^{e_i})})$.

**Proof.** In addition to defining the length of the period of $\{PN_n^{(p_i^{e_i})}\}$ with $\text{Per}(PN^{(p_i^{e_i})})$, we also express this statement the sequence $\{PN_n^{(p_i^{e_i})}\}$ repeats only after blocks of length

$u. \mathrm{Per}(PN^{(p_i^{e_i})}), u \in \mathbb{N}$ and this definition is due to the fact that we have already shown the length of the period of $\{PN_n^{(m)}\}$ with $\mathrm{Per}(PN^{(m)})$. This means that for all values of $i$, $\{PN_n^{(p_i^{e_i})}\}$ repeats after $\mathrm{Per}(PN^{(m)})$.

Thus, $\mathrm{Per}(PN^{(m)})$ is of the form $u.\mathrm{Per}(PN^{(p_i^{e_i})})$ for all values of $i$, and any such number gives a period of $\{PN_n^{(m)}\}$. Then we get that $\mathrm{Per}(PN^{(m)}) = \mathrm{lcm}[\mathrm{Per}(PN^{(p_i^{e_i})})]$, as required.

**Conclusion 4.** For $p_i$, if it consists of even numbers such as 2, 4 or 6, simple periodicity is provided, otherwise periodicity is provided. We can illustrate this with examples below.

**Example 5.** $\{PN_n^{(2)}\} = \{0,1,1,0,1,1,...\} \Rightarrow \mathrm{Per}(PN^{(2)}) = 3$.

For $m = 4, e_1 = 2, \{PN_n^{(4)}\} = \{0,1,1,2,1,3,0,1,1,2,1,3,0,1,1,...\} \Rightarrow \mathrm{Per}(PN^{(4)}) = 6$

For $m = 3, \{PN_n^{(3)}\} =$
$\{0,1,1,2,2,2,0,2,0,0,2,1,2,0,1,1,1,0,1,0,0,1,2,1,0,2,2,1,1,1,0,1,0,0,1,2,1,...\} \Rightarrow$
$\mathrm{Per}(PN^{(3)}) = 13$

For $m = 6, \{PN_n^{(6)}\} =$
$\{0,1,1,2,5,5,0,5,3,0,5,1,2,3,1,4,5,5,2,3,5,0,3,5,4,5,3,4,1,5,2,5,3,2,3,3,2,1,5,0,1,1,2,5,5,0,...\} \Rightarrow$
$\mathrm{Per}(PN^{(6)}) = 39$

$\mathrm{Per}(PN^{(6)}) = \mathrm{lcm}\left(\mathrm{Per}(PN^{(2)}), \mathrm{Per}(PN^{(3)})\right) = \mathrm{lcm}(3,13) = 39$

as required.

For $m = 5, \{PN_n^{(5)}\} =$
$\{0,1,1,2,0,1,4,3,2,3,4,0,3,0,0,3,2,4,1,4,2,0,40,0,4,3,1,1,0,1,3,1,3,4,4,1,1,1,3,2,0,3,3,1,0,3,2,4,1,...\} \Rightarrow$
$\mathrm{Per}(PN^{(5)}) = 31$

For the matrix $\boldsymbol{A} = [a_{ij}]_{(k+1)\times(k+1)}$ with $a_{ij}$ integers, $A(mod\ m)$ means that all entries of $A$ are reduced modulo $m$, that is, $A(mod\ m) = (a_{ij}(mod\ m))$. Let $\langle PN \rangle_{p^\alpha} = \{K^i (mod\ p^\alpha) | i \geq 0\}$ be a cyclic group and $|\langle PN \rangle_{p^\alpha}|$ denote the order of $\langle PN \rangle_{p^\alpha}$. From (5), we have that $\mathrm{Per}(PN^{(p^\alpha)}) = |\langle PN \rangle_{p^\alpha}|$.

**Theorem 6.** Let $t$ be the positive integer such that $\text{Per}(PN^{(p)}) = \text{Per}(PN^{(p^t)})$. Then we have $\text{Per}(PN^{(p^\alpha)}) = p^{\alpha-t}\text{Per}(PN^{(p)})$, $\alpha \geq t$. In particular, if $\text{Per}(PN^{(p)}) \neq \text{Per}(PN^{(p^2)})$ then we have $\text{Per}(PN^{(p^\alpha)}) = p^{\alpha-1}\text{Per}(PN^{(p)})$, $\alpha > 1$.

**Proof.** Let $q$ be a positive integer. Since $K^{\text{Per}(N^{(p^{q+1})})} \equiv I(mod\ p^{q+1})$ and $K^{\text{Per}(N^{(p^{q+1})})} \equiv I(mod\ p^q)$, we get that $\text{Per}(PN^{(p^q)})$ divides $\text{Per}(PN^{(p^{q+1})})$ where $I$ is an identity matrix. On the other hand, we know

$$K^{\text{Per}(PN^{(p^q)})} = I + (a_{ij}^{(q)}p^q).$$

So, we have

$$K^{\text{Per}(PN^{(p^q)})p} = (I + a_{ij}^{(q)}p^q)^p = \sum_{i=0}^{p}\binom{p}{i}(a_{ij}^{(q)}p^q)^i \equiv I(mod\ p^{q+1}),$$

which yields the result that $\text{Per}(PN^{(p^{q+1})})$ divides $\text{Per}(PN^{(p^q)}p)$. Therefore, we get

$$\text{Per}(PN^{(p^{q+1})}) = \text{Per}(PN^{(p^q)}) \text{ or } \text{Per}(PN^{(p^{q+1})}) = \text{Per}(PN^{(p^q)}p)$$

and $\text{Per}(PN^{(p^{q+1})}) = \text{Per}(PN^{(p^q)}p)$ if and only if there is an $a_{ij}^{(q)}$ which is not divisible by $p$. Since $\text{Per}(PN^{(p^t)}) \neq \text{Per}(PN^{(p^{t+1})})$, there is an $a_{ij}^{(t+1)}$ which is not divisible by $p$. Thus, we get $\text{Per}(PN^{(p^{t+1})}) \neq \text{Per}(PN^{(p^{t+2})})$. Then by an inductively on $t$ desired is achieved.

Let us explain this theorem with an example.

**Example 7.** For $p = 2$ and $q = 1$, $K^{\text{Per}(PN^{(4)})} \equiv K^6 \equiv I(mod\ 4)$,

so that,

$$K^6 = \begin{pmatrix} 9 & 4 & 20 \\ 20 & 9 & 44 \\ 44 & 20 & 97 \end{pmatrix}_{mod 4} = \begin{pmatrix} 1 & 0 & 0 \\ 0 & 1 & 0 \\ 0 & 0 & 1 \end{pmatrix} = I.$$

Also

$$K^{\text{Per}(N^{(2)})2} = \left(\begin{pmatrix} 1 & 0 & 0 \\ 0 & 1 & 0 \\ 0 & 0 & 1 \end{pmatrix} + 2\begin{pmatrix} 4 & 2 & 10 \\ 10 & 4 & 22 \\ 22 & 10 & 48 \end{pmatrix}\right)^2 = \begin{pmatrix} 1 & 0 & 0 \\ 0 & 1 & 0 \\ 0 & 0 & 1 \end{pmatrix}(mod\ 2^2).$$

In this example, for $t = 1$ and $\alpha = 2$, we have

$$\text{Per}(PN^{(p^\alpha)}) = p^{\alpha-t}\text{Per}(PN^{(p)}) = \text{Per}(PN^{(2^2)}) = 2^1\text{Per}(PN^{(2)}),$$

where $\text{Per}(PN^{(2)}) = 3 \ and \ \text{Per}(PN^{(2^2)}) = 6$ from Example 7. Also, in this example for $\text{Per}(PN^{(p)}) \neq \text{Per}(PN^{(p^2)})$, $\text{Per}(PN^{(p^\alpha)}) = p^{\alpha-1}\text{Per}(PN^{(p)})$ is provided.

## 3. Pell-Narayana orbits of finite groups

**Definition 8.** For $j$-generating $(2 \leq j \leq 3)$, the Pell-Narayana orbits of finite groups are defined as follows:

i) Let $(x_0, x_1) \in X$ be the generating pair for a 2-generator group $G$, in that case, the sequence $\{s_i\}$ of elements of $G$ defines the Pell-Narayana orbit $PN(G)_{(x_0,x_1)}$:

$$s_0 = x_0, s_1 = x_1, s_2 = (s_1)^2, \dots, s_n = (s_{n-3})(s_{n-1})^2 \text{ for } n \geq 3.$$

ii) Let $(x_0, x_1, x_2) \in X$ be the generating triplet for a 3-generator group $G$, in that case, the sequence $\{s_i\}$ of elements of $G$ defines the Pell-Narayana orbit $PN(G)_{(x_0,x_1,x_2)}$:

$$s_0 = x_0, s_1 = x_1, s_2 = x_2, \dots, s_n = (s_{n-3})(s_{n-1})^2 \text{ for } n \geq 3.$$

**Theorem 9.** A Pell-Narayana orbit of a finite group with $j$-generating $(2 \leq j \leq 3)$ is simply periodic sequence.

**Proof.** Let's consider group $G$ with the generating pair $(x_0, x_1)$. If the order of the group $G$ is $n$, we can say that there are $n^3$ different 3-tuples elements in $G$.

This means that every 3-tuples element in the Pell-Narayana orbit of $G$ will appear at least 2 times.

This means very clearly that the Pell-Narayana orbit is periodic sequence.

Given the periodicity here, there are positive integer $i$ and $j$ with $i > j$ such that

$$s_{i-1} = s_{j-1}, s_{i-2} = s_{j-2}, s_{i-3} = s_{j-3}.$$

From the definition of the orbit in the Pell-Narayana sequence, we know that

$$s_i = (s_{i-3})(s_{i-1})^2 \text{ and } s_j = (s_{j-3})(s_{j-1})^2.$$

Thus, $s_i = s_j$ and it follows that

$$s_{i-j} = s_{j-j} = s_0 = x_0,$$

$$s_{i-j+1} = s_{j-j+1} = s_1 = x_1.$$

So, the Pell-Narayana orbit is simply periodic sequence. The proof of 3-generator groups is similarly done.

We denote the periods of the orbits $PN(G)_{(x_0,\ldots,x_n)}$ with $1 \leq n \leq 2$ by $LPN(G)_{(x_0,\ldots,x_n)}$.

Then it is clear from the operations that the Pell-Narayana orbits of the finite groups depend on the selected generating set and the order for the assignments of $x_0, \ldots, x_n; 1 \leq n \leq 2$.

**Definition 10.** For a finite group $G$, group $G$ is a Pell-Narayana sequenceable if $G$ has a Pell-Narayana orbit such that every element of the $G$ group appears in the sequence.

Based on the definition, it is said that all the examples given below are a Pell-Narayana sequenceable as a requirement of the expression given in Definition 10.

We examine Pell-Narayana orbits of finite groups, $Q_8, Q_8 \times \mathbb{Z}_{2m}$ and $Q_8 \times_\varphi \mathbb{Z}_{2m}$, $m \geq 3$.

**Theorem 11.** $LPN(Q_8)_{(x,y)} = 6$.

**Proof.** $LPN(Q_8)_{(x,y)}$ is

$$x, y, y^2, x, y^3, e, x, y, y^2, x, y^3, \ldots$$

which has period $LPN(Q_8)_{(x,y)} = 6$.

**Theorem 12.** For for each generating triplet, the period of the Pell-Narayana orbit of the group $Q_8 \times \mathbb{Z}_{2m}, (m \geq 3)$ is $lcm[6, \text{Per}(PN^{(2m)})]$.

**Proof.** Consider the Pell-Narayana orbit $PN(Q_8 \times \mathbb{Z}_{2m})_{(x,y,z)}$:

$$x, y, z^{a_1}, xz^{a_2}, y^3z^{a_3}, y^2z^{a_4}, xz^{a_5}, yz^{a_6}, z^{a_7}, xz^{a_8}, y^3z^{a_9}, y^2z^{a_{10}}, xz^{a_{11}}, yz^{a_{12}}, z^{a_{13}}, xz^{a_{14}}, \ldots$$

Using the above information, we have

$s_0 = x, s_1 = y, s_2 = z,$

$s_3 = xz^2, s_4 = y^3z^3, s_5 = y^2z^4, s_6 = xz^5, s_7 = yz^6, s_8 = z^7,$

$s_9 = xz^8, s_{10} = y^3z^9, s_{11} = y^2z^{10}, s_{12} = xz^{11}, s_{13} = yz^{12}, s_{14} = z^{13},$

$s_{6i-3} = xz^{6i-4}, s_{6i-2} = y^3z^{6i-3}, s_{6i-1} = y^2z^{6i-2}, s_{6i} = xz^{6i-1}, s_{6i+1} = yz^{6i},$

$s_{6i+2} = z^{6i+1}, \ldots$

The sequence has the form layers of length 6. So, we need an $i$ such that $s_{6i} = x$, $s_{6i+1} = y$, $s_{6i+2} = z$. It is seen that the Pell-Narayana orbit $PN(Q_8 \times \mathbb{Z}_{2m})_{(x,y,z)}$ has period $lcm[6, \text{Per}(PN^{(2m)})]$.

**Theorem 13.** For each generating triplet, the period of the Pell-Narayana orbit of the group $Q_8 \times_\varphi \mathbb{Z}_{2m}$, $(m \geq 3)$ is $lcm[6, \text{Per}(PN^{(2m)})]$.

**Proof.** Consider the Pell-Narayana orbit $PN(Q_8 \times_\varphi \mathbb{Z}_{2m})_{(x,y,z)}$:

$$x, y, z^{a_1}, xz^{a_2}, y^3 z^{a_3}, y^2 z^{a_4}, xz^{a_5}, yz^{a_6}, z^{a_7}, xz^{a_8}, y^3 z^{a_9}, y^2 z^{a_{10}}, xz^{a_{11}}, yz^{a_{12}}, z^{a_{13}}, xz^{a_{14}}, \ldots$$

From here, we get the following sequence.

$s_0 = x, s_1 = y, s_2 = z$,

$s_3 = xz^2, s_4 = y^3 z^3, s_5 = y^2 z^4, s_6 = xz^5, s_7 = yz^6, s_8 = z^7$,

$s_9 = xz^8, s_{10} = y^3 z^9, s_{11} = y^2 z^{10}, s_{12} = xz^{11}, s_{13} = yz^{12}, s_{14} = z^{13}$,

$s_{6i-3} = xz^{6i-4}, s_{6i-2} = y^3 z^{6i-3}, s_{6i-1} = y^2 z^{6i-2}, s_{6i} = xz^{6i-1}, s_{6i+1} = yz^{6i}$,

$s_{6i+2} = z^{6i+1}, \ldots$

The sequence has the form layers of length 6. So, we need an $i$ such that $s_{6i} = x$, $s_{6i+1} = y$, $s_{6i+2} = z$. So, we can see that the Pell-Narayana orbit $PN(Q_8 \times_\varphi \mathbb{Z}_{2m})_{(x,y,z)}$ has period $lcm[6, \text{Per}(PN^{(2m)})]$.

The other generating triplets are proved in a similar way.

**Remarks 14.** If $\text{Per}(PN^{(2m)}) \leq 2m$ and $6 | \text{Per}(PN^{(2m)})$, the groups $Q_8 \times \mathbb{Z}_{2m}$ and $Q_8 \times_\varphi \mathbb{Z}_{2m}$ such that $m \geq 3$ are not Pell-Narayana sequence.

## 4. The Pell-Narayana Length of Generating Pairs in Groups

**Theorem 15.** The Pell-Narayana length of $\mathbb{Z}_n \times \mathbb{Z}_m$ (where $\mathbb{Z}_n = \langle x \rangle$ and $\mathbb{Z}_m = \langle y \rangle$) equals $lcm[\text{Per}(PN^{(n)}), \text{Per}(PN^{(m)})]$.

**Proof.** We know that $\mathbb{Z}_n \times \mathbb{Z}_m$ have the following presentation

$$\langle x, y : x^n = y^m = e, xy = yx \rangle.$$

The Pell-Narayana orbit is:

$$s_0 = x, s_1 = y, s_2 = y^2, s_3 = xy^4, s_4 = x^2y^9, s_5 = x^4y^{20}, s_6 = x^9y^{44}, \ldots$$

If we get, $x_i = x, x_{i+1} = y, x_{i+2} = y^2$ then the proof is finished.

Examining this statement in more detail, it yields

$$x^{PN_{i-2}}y^{PN_i} = x = s_0, \; x^{PN_{i-1}}y^{PN_{i+1}} = y = s_1, \; x^{PN_i}y^{PN_{i+2}} = y^2 = s_2.$$

The least non-trivial integer satisfying the above conditions occurs when

$$i = lcm[\text{Per}(PN^{(n)}), \text{Per}(PN^{(m)})].$$

**4.1 Applications**

In this section, let $\ell, m, n > 1$ be integers.

**Definition 16.** The polyhedral group $(\ell, m, n)$ is defined by

$$\langle x, y, z : x^\ell = y^m = z^n = xyz = e \rangle \text{ or } \langle x, y : x^\ell = y^m = (xy)^n = e \rangle.$$

The polyhedral group $(\ell, m, n)$ is finite if and only if $k = \ell mn \left(\frac{1}{\ell} + \frac{1}{m} + \frac{1}{n} - 1\right) = mn + n\ell + \ell m - \ell mn$ is positive. Its order is $\frac{2\ell mn}{k}$.

**Definition 17.** The binary polyhedral group $\langle \ell, m, n \rangle$ is given by

$$\langle x, y, z : x^\ell = y^m = z^n = xyz \rangle, \text{ or } \langle x, y : x^\ell = y^m = (xy)^n \rangle.$$

The binary polyhedral group $\langle \ell, m, n \rangle$ is finite if and only if the number

$$k = \ell mn \left(\frac{1}{\ell} + \frac{1}{m} + \frac{1}{n} - 1\right) = mn + n\ell + \ell m - \ell mn,$$

is positive. Its order is $\frac{4\ell mn}{k}$.

We find the Narayana lengths of the polyhedral groups $(2,2,2), (n, 2,2), (2, n, 2), (2,2, n)$

and the binary polyhedral groups $\langle 2,2,2 \rangle, \langle n, 2,2 \rangle, \langle 2, n, 2 \rangle, \langle 2,2, n \rangle$ for the generating pair $(x, y)$.

**Theorem 18.** The Pell-Narayana length of the polyhedral group $(2,2,2)$ is 3.

**Proof.** From Theorem 15, we can see that $LPN_{x,y,y}((2,2,2)) = 3$. Since $(2,2,2) \cong \mathbb{Z}_2 \times \mathbb{Z}_2$, we have

$$s_0 = x, s_1 = y, s_2 = e, s_3 = x, s_4 = y, s_5 = e, \ldots$$

and $LPN_{x,y,y}((2,2,2)) = 3$.

**Theorem 19.** The Pell-Narayana length of the polyhedral group $(2, n, 2)$ is

$$LPN_{x,y,y}((2,n,2)) = \begin{cases} \frac{3n}{2}, & n \equiv 0 \pmod{2} \\ 3n, & n \equiv 1 \pmod{2} \end{cases}, n > 2$$

**Proof.** If $\langle x, y: x^2 = y^n = (xy)^2 = e \rangle$, $|x| = 2$, then $|y| = n$ and $|xy| = 2$.

The Pell-Narayana orbit is:

$$x, y, y^2, xy^4, y, y^4, xy^{12}, y, y^6, xy^{24}, y, y^8, xy^{40}, \ldots$$

$s_0 = x, s_1 = y, s_2 = y^2, s_3 = xy^4, s_4 = y, s_5 = y^4, s_6 = xy^{12}, s_7 = y, s_8 = y^6,$

$s_9 = xy^{24}, s_{10} = y, s_{11} = y^8, \ldots, s_{3i} = xy^{2(i^2+i)}, s_{3i+1} = y, s_{3i+2} = y^{2(i+1)}$.

so, we get $2i = un$ for $u, i \in \mathbb{N}$.

If $n \equiv 0 \pmod{2}$ then $n = 2k$ and $2i = u2k$, $k = \frac{i}{u}$, $n = \frac{2i}{u}$, $nu = 2i$, $i = \frac{nu}{2}$ for $u = 1, i = \frac{n}{2}$.

So, we get

$$LPN_{x,y,y}((2,n,2)) = lcm\left(\left(3, \frac{n}{2}, 3\right)\right) = 3\frac{n}{2}.$$

If $n \equiv 1 \pmod{2}$ then $n = 2k+1$ and $2i = u(2k+1)$, $k = \frac{i-1}{2}$, $n = \frac{2(i-1)}{2} + 1$, $n = i$.

So, we get

$$LPN_{x,y,y}((2,n,2)) = lcm(3, n, 3) = 3n.$$

**Theorem 20.** Let $G$ be anyone of the polyhedral group $(n, 2,2)$ and $(2,2, n)$. Then

$$LPN_{x,y,y}(G) = \begin{cases} \frac{3n}{2}, & n \equiv 0 \pmod{2} \\ 3n, & n \equiv 1 \pmod{2} \end{cases}, n > 2.$$

**Proof.** Let us consider the group $(n, 2,2)$. We know that

$\langle x, y: x^n = y^2 = (xy)^2 = e \rangle$, $|x| = n$, $|y| = 2$ and $|xy| = 2$. The Pell-Narayana orbit is

$s_0 = x, s_1 = y, s_2 = e, s_3 = x, s_4 = yx^2, s_5 = e, s_6 = x, s_7 = yx^4, s_8 = e, \ldots, s_{3i} = x, s_{3i+1} = yx^{2i}, s_{3i+2} = e.$

So, we want to find $i \in \mathbb{N}$ such that $2i = un$ for $u \in \mathbb{N}$.

If $n \equiv 0 (mod 2)$ then $n = 2k$ and $2i = u2k$, $k = \frac{i}{u}$

$n = \frac{2i}{u}, nu = 2i, i = \frac{nu}{2}$ for $u = 1, i = \frac{n}{2}$.

So, we get

$$LPN_{x,y,y}((n,2,2)) = lcm\left(\left(\frac{n}{2}, 3, 3\right)\right) = 3\frac{n}{2}.$$

If $n \equiv 1 (mod 2)$ then $n = 2k + 1$ and $2i = u(2k+1)$, $k = \frac{i-1}{2}$

$n = \frac{2(i-1)}{2} + 1, n = i$.

So, we get

$$LPN_{x,y,y}((n,2,2)) = lcm(n,3,3) = 3n.$$

**Example 21.** For $n = 4$, the Pell-Narayana length of the polyhedral group $(4,2,2)$ is 6. So that

$$x, y, e, x, yx^2, e, \ldots$$

**Theorem 22.** The Pell-Narayana length the binary polyhedral group $\langle 2,2,2 \rangle$ is 6.

**Proof.** Since the group has the presentation $\langle x, y: x^2 = y^2 = (xy)^2 \rangle$ and $|x| = 4$, the Pell-Narayana orbit is $, y, y^2, x, y^3, e, x, y, y^2, x, y^3, e, \ldots$. So, we get $LN_{x,y,y}((2,2,2)) = 6$.

**Theorem 23.** For $n \geq 2$, the Pell-Narayana length of the binary polyhedral group $\langle 2, n, 2 \rangle$ is

$$LPN_{x,y,y}(\langle 2, n, 2\rangle) = \begin{cases} 3n, & n \equiv 0 (mod\ 2) \\ 6n, & n \equiv 1 (mod\ 2) \end{cases}.$$

**Proof.** Since $\langle x, y: x^2 = y^n = (xy)^2 \rangle$, $|x| = 4, |y| = 2n, |xy| = 4$, the Pell-Narayana orbit is

$$s_0 = x, s_1 = y, s_2 = y^2, s_3 = xy^4, s_4 = yx^2, s_5 = y^4, s_6 = xy^{12}, s_7 = y, s_8 = y^6, s_9 = xy^{24}, s_{10} = yx^2, s_{11} = y^8, \ldots$$

For $n \equiv 0 (mod\ 2)$,

$$s_{3i-1} = y^{2i}, s_{3i} = xy^{2(i^2+i)}, s_{3i+1} = y,$$

and for $n \equiv 1 (mod\ 2)$,

$$s_{3i-1} = y^{2i}, s_{3i} = xy^{2(i^2+i)}, s_{3i+1} = yx^2.$$

Now we want to find $i \in \mathbb{N}$ such that $2i = un\ for\ u \in \mathbb{N}$.

For first case:

$$n = 2k, i = uk.$$

For $u = 2, i = n$, we have

$$LPN_{x,y,y}(\langle 2, n, 2 \rangle) = lcm(3, n, 3) = 3n.$$

For second case:

$$n = 2k + 1, 2i = u(2k + 1).$$

For $u = 4, i = 2n$, we get

$$LPN_{x,y,y}(\langle 2, n, 2 \rangle) = lcm(3, 2n, 3) = 6n.$$

**Theorem 24.** Let $G$ be any of these two groups, $\langle n, 2, 2 \rangle$ and $\langle 2, 2, n \rangle$ then

$$LPN_{x,y,y}(G) = \begin{cases} 3n, & n \equiv 0 (mod\ 2) \\ 6n, & n \equiv 1 (mod\ 2) \end{cases}, n > 2.$$

**Proof.** Now, let $G$ be the binary polyhedral group $\langle n, 2, 2 \rangle$. The group is defined by

$$\langle x, y : x^n = y^2 = (xy)^2 \rangle, |x| = 2n, |y| = 4, |xy| = 4.$$

The Pell-Narayana orbit is:

$$s_0 = x, s_1 = y, s_2 = y^2, s_3 = x, s_4 = yx^2, s_5 = e, s_6 = x, s_7 = yx^4, s_8 = y^2, s_9 = x, s_{10} = yx^6, s_{11} = e, \ldots$$

For $n \equiv 0 (mod\ 2)$, we get

$$s_{3i} = x, s_{3i+1} = xy^{2i}, s_{3i+2} = e.$$

So, we need $i \in \mathbb{N}$ such that $2i = un$ for $u \in \mathbb{N}$.

For first case:

$$n = 2k, i = uk.$$

For $u = 2, i = n$, we get

$$LPN_{x,y,y}(\langle n, 2, 2 \rangle) = lcm(n, 3, 3) = 3n.$$

For second case:

$$n = 2k + 1, 2i = u(2k + 1).$$

For $u = 4, i = 2n$, we get

$$LPN_{x,y,y}(\langle n, 2, 2 \rangle) = lcm(2n, 3, 3) = 6n.$$

Now, let G be the binary polyhedral group $\langle 2,2,n \rangle$. We first note that in the group defined by

$$\langle x, y: x^2 = y^2 = (xy)^n \rangle, |x| = 4, |y| = 4, |xy| = 2n.$$

The proofs are like that of firstly case.

**Example 25.**

In the case $n = 4$, the Pell-Narayana length of the binary polyhedral group $\langle 4,2,2 \rangle$ is 12.

The group is defined by

$$\langle x, y: x^4 = y^2 = (xy)^2 \rangle, |x| = 8, |y| = 4, |xy| = 4.$$

So, the Pell-Narayana sequence in the $\langle 4,2,2 \rangle$ is

$$x, y, y^2, x, yx^2, e, x, yx^4, y^2, x, yx^6, e, x, y, y^2, x, yx^2, \ldots$$

## 5. CONCLUSIONS

We have recalled the essential features of the Pell-Narayana sequence and the main properties so that we could examine the Narayana sequences modulo *m*. So, we have defined the Pell-Narayana orbits of 2-generator and 3-generator finite groups. Also, the period of the Pell-Narayana orbit of the groups $Q_8, Q_8 \times \mathbb{Z}_{2m}$ and $Q_8 \times_\varphi \mathbb{Z}_{2m}$ for $m \geq 3$ were obtained. Finally, we have obtained the Pell-Narayana lengths of the polyhedral and the binary polyhedral groups with specific examples.